# Category theoretic analysis of hierarchical protein materials and social networks

*David I. Spivak[1], Tristan Giesa[2,3], Elizabeth Wood[2] and Markus J. Buehler[2,4,5]\**

[1] *Department of Mathematics, Massachusetts Institute of Technology, 77 Massachusetts Ave. Cambridge, MA, USA*

[2] *Laboratory for Atomistic and Molecular Mechanics, Department of Civil and Environmental Engineering, Massachusetts Institute of Technology, 77 Massachusetts Ave. Room 1-235A&B, Cambridge, MA, USA*

[3] *Department of Mechanical Engineering, RWTH Aachen University, Templergraben 55, 52056 Aachen, Germany*

[4] *Center for Materials Science and Engineering, Massachusetts Institute of Technology, 77 Massachusetts Ave., Cambridge, MA, USA*

[5] *Center for Computational Engineering, Massachusetts Institute of Technology, 77 Massachusetts Ave., Cambridge, MA, USA*

*\* Corresponding author: [mbuehler@MIT.EDU](mailto:mbuehler@MIT.EDU)*

**Abstract**: Materials in biology span all the scales from Angstroms to meters and typically consist of complex hierarchical assemblies of simple building blocks. Here we describe an application of category theory to describe structural and resulting functional properties of biological protein materials by developing so-called ologs. An olog is like a "concept web" or "semantic network" except that it follows a rigorous mathematical formulation based on category theory. This key difference ensures that an olog is unambiguous, highly adaptable to evolution and change, and suitable for sharing concepts with other olog. We consider simple cases of alpha-helical and amyloid-like protein filaments subjected to axial extension and develop an olog representation of their structural and resulting mechanical properties. We also construct a representation of a social network in which people send text-messages to their nearest neighbors and act as a team to perform a task. We show that the olog for the protein and the olog for the social network feature identical category-theoretic representations, and we proceed to precisely explicate the analogy or isomorphism between them. The examples presented here demonstrate that the intrinsic nature of a complex system, which in particular includes a precise relationship between structure and function at different hierarchical levels, can be effectively represented by an olog. This, in turn, allows for comparative studies between disparate materials or fields of application, and results in novel approaches to derive functionality in the design of *de novo* hierarchical systems. We discuss opportunities and challenges associated with the description of complex biological materials by using ologs as a powerful tool for analysis and design in the context of materiomics, and we present the potential impact of this approach for engineering, life sciences, and medicine.

## 1. Introduction

Biological materials, many of which contain proteins as a basic building block, provide an enormous diversity of properties including structural support, prey procurement and material transport [1]. Significant evidence has now emerged that proteins are organized in functional networks, resulting in structures that span many hierarchical scales [2,3,4,5,6,7,8,9]. In the glassy sponge *Eucleptella*, for instance, silica nanospheres are arranged at multiple levels of hierarchy to constitute a skeleton with high structural stability at minimum cost [10]. The teeth of sea urchins and the lamellar structure of mollusk shells are other examples for structural hierarchies in biomaterials that lead to extremely strong and tough structures [11]. Earlier studies showed that in materials like bone or wood, for example, the structural assembly of basic building blocks such as collagen, water, hydroxyapatite minerals, hemicelluloses and lignin governs the mechanical properties at different length scales with similar mechanisms despite the differences in the building blocks and the overall material properties [12,13]. A frontier in protein materials science is the understanding how the exceptionally complex functionality found in natural biological systems is created despite *i)* a limited number of 20 amino acid building blocks, *ii)* constraints in available material volume and energy for synthesis, and *iii)* only a handful of simple chemical interaction force fields, generally referred to as interaction rules [2,3,4,5,6,14].

It is remarkable that the same library of amino acid building blocks creates materials as diverse as spider silk, tendon, cornea, blood vessels, and cellular protein networks, each of which displays greatly variegated functions. Our understanding of the synthesis of their basic elements into multi-functional structures remains in its infancy, and is currently limited to specific protein networks or protein materials. For example, mechanistic theories are typically developed for specific proteins (see *e.g.* [15,16,17,18,19,20]) rather than providing a unified model that is applicable to a variety of distinct materials. The extraction of generic principles of how functional properties are derived in functionally diverse systems despite the presence of the same (universal) building blocks, solely by using structure as a design paradigm, presents an exciting opportunity. The systematic characterization of this knowledge has resulted in the formation of a new field referred to as materiomics [21].

Here we describe, by means of the application of the mathematical field of category theory to protein materials, how the extreme diversity of protein functional properties can be described in a unified model that contains only a limited number of universal elements and their interaction rules. Category theory has been successfully applied to carry out a qualitative analysis in fields such as linguistics (grammar, syntax, semantics, *etc.* – key concepts that enable the understanding of language, see *e.g.* [22,23,24,25,26,27]) and computer science (again modeling syntax and semantics of denotation and operation in programming languages, see *e.g.* references [28,29,30,31]).

Category theory can be seen as an abstraction of graph theory which has been used to describe the structure of biopolymers, disease spread and neuronal activity as well as to determine the role of proteins or genes of unknown function and to identify drug targets [32,33,34,35,36,37,38,39]. The focus of such earlier studies has been on metabolic networks of biochemical reactions, protein interaction networks and transcriptional regulatory networks, amongst others [32]. The role of network motifs (or building blocks) has been studied but is often limited to structural motifs in networks [40]. Current theoretical approaches to material science focus either on structural aspects or on functional aspects and lack the general and abstract description of how system elements behave and interact with each other in order to create functionality. For protein materials, a major limitation of existing methods is the reliance on topological considerations which do not account for biological information about the network's interaction [32]. For example, protein networks are typically modeled as undirected graphs where the nodes represent proteins and the edges represent physical interactions between them. The need for a broader view on functionality (biological, chemical, mechanical and other) and structure, enhanced data collection abilities, and integrated studies, as each study covers only a small subset of a generally big topic, cannot be achieved with graph theoretic approaches. Category theory provides means to overcome the limitations of conventional networks while also including graph theoretic tools.

The combination of universal elements into multi-level structures enables protein materials to achieve context specific functionalities in an abstract "complexity space". This paradigm showed that in order to create highly functional materials it is not essential to rely on a multitude and a certain quality of building blocks (*e.g.* with superior qualities, great material volume, strong interactions, *etc.*). Rather, it is sufficient to use simple interaction rules and simple building blocks – each of which does not need to possess superior qualities – but assembled into hierarchical systems, where the overall structure provides enhanced functionality [14,21]. This insight has implications for our understanding of how nanomaterials could be utilized to create macroscopically functional materials, and suggests a paradigm that departs greatly from the current belief in engineering science that material building blocks with superior qualities at the small scale (*e.g.*, carbon nanotube, carbyne, graphene, *etc.*) are critical to reach high performance materials. On the contrary, we hypothesize that superior functionality can be reached with any fundamental building block, provided that the design space is expanded to incorporate hierarchical structures. Eventually, an understanding of how diverse functional properties can arise out of inferior building blocks could make a profound impact towards the development of environmentally benign and friendly materials, as it would allow manufacturers to use local, abundant, and simple building blocks with overall negative $CO_2$ balance (*e.g.* wood, plants, silica, water, soy beans) to create highly functional materials and structures. But how can we find a proper mathematical description of these hierarchical mechanisms that generate functional properties? A possible approach is to use novel mathematical concepts that provide a powerful, abstract way to describe emergence of functionality from first principles, *e.g.* on the basis of fundamental interactions between building blocks.

**1.1 Abstract representation of structure and function of protein materials using category theory**

Biological materials evolved to perform specific biological functions [2,3,4,5,6,7,8,9], where the components and connections within a given biological material are analogous to a circuit diagram. But just as it is extremely difficult to deduce the circuitry of a device by experimenting with its inputs and outputs, it is similarly inadequate to describe the higher-level structure of a biological material using only the physical interactions between proteins and some information about gene expression. Instead, we need to take into account additional types of structural information given by the fundamental principles that govern the interactions of the building blocks that define the system and its emerging functionality as these building blocks are connected together, from the micro to the macro scale. The above considerations are important in any synthetic science; in order to duplicate the functionality of a natural system, we do not need to understand everything about it, only the principles out of which the desired functions arise. Biological systems contain any number of copies of thousands of different components, each with very specific interactions, and each representing a microscopic device in and of itself. As a result, the microscopic description of a biological system (and materials therein) is intractably complex, unless one moves to a higher level of abstraction in the analysis that, as discussed before, cannot be solely provided by network theory.

It is exactly in the face of this complexity that ologs are so appealing. The olog presents us with an opportunity to identify patterns that describe systems and their components, to elucidate possible connections among these components, and to construct isolated *functional* (and specifically not limited to structural) "modules" by comparing information from many different materials or organisms. That is, by determining fundamental design principles that are simple yet functional, we can not only produce a powerful conceptual model of our system, we also create the possibility of comparing vastly different systems. Indeed, we will show below that although there is almost nothing physical in common between a protein and a social network, we can construct a scenario in which the design principles are well-matched, and thus the systems may be compared. Such a comparison may allow results from the science of social networks to guide us in our study of biological materials of the same structure, and vice versa.

To give a few concrete examples of how such analogies between seemingly disparate fields can be made, Figure 1 shows an illustration of multiscale hierarchical structure of protein materials, a summary of multiscale modeling and experimental tools, and an analogy to music. In protein materials (left for the example of spider silk), multifunctional materials are defined *via* the formation of hierarchical structures. The synergistic interaction of structures and mechanisms at multiple scales provides the basis for enhanced functionality of biological materials despite the reliance on few distinct building blocks. Similar to the case of protein materials, musical composition (right) is built upon universal elements at the microscale such as basic wave forms, and gathers a small variety of available instruments into hierarchical assemblies to create macroscale functionality, such as a particular orchestral sound (*e.g.* a symphony). Universality tends to dominate at smaller levels, whereas diversity is found predominantly at larger, functional levels [9]. The integrated use of computational and experimental methods at multiple scales provides a powerful approach to elucidate the scientific concepts underlying the materiomics paradigm (center).

**1.2 Outline of this paper**

The scope of this paper is to present a novel methodology to material science which incorporates structural and functional hierarchies. Hence, we utilize a comprehensive example, the behavior of an amyloid and an alpha-helical structure under tensile load, in order to illustrate the concept rather than to derive a full analysis. The biochemical structure is extremely simplified, thus allowing us to demonstrate the transformation of the protein system into a social network with similar characteristics. We conclude with a discussion of opportunities for the science and engineering of natural protein materials as well as synthetically designed materials from the atomistic scale with the chemical structure of molecules to the macroscale.

**2. Methods**

Category theory is a relatively new branch of mathematics (invented 200 years after the introduction of partial differential equations), designed to connect disparate fields within the larger discipline (see [41]). It is both a language that captures the essential features of a given subject, and a toolbox of theorems that can be applied quite generally. If a given study within mathematics is formalized as a category, it can be connected with other categories that are seemingly far afield, as long as these structures align in the required "functorial" way. Theorems within one branch, like abstract equational algebra, can be applied to a totally different area, like geometric topology. Category theory may not only serve as an alternate foundation to mathematics [42], it unites the various distinct areas within advanced mathematics, formally proving similarities that are not apparent on the surface [43]. A good overview for non-specialists can be found in [44] and [45].

Quickly after its inception, category theorists realized that its basic ideas were applicable well beyond the borders of mathematics. Category theory has by now been successfully applied in computer science, linguistics, and physics [46]. Whereas the theory of differential equations can be applied throughout science to create *quantitative models*, category theory can be applied throughout science to create *qualitative models*. And once such a qualitative model is formed as a category, its basic structure can be meaningfully compared (again *via functors*) with that of any other category, be it mathematical, linguistic, or other [47]. Like a biological system, the basic building blocks of a category are simple, but the networks that can be formed out of them are as complex as mathematics itself. These building blocks are called objects, arrows, and composition: arrows between objects form paths which can be composed into new arrows. It is a wonder that such a simple system can account for the wide variety of forms found in the mathematical universe, but perhaps this is less of a surprise to a biologist who notices the same phenomenon in his or her field.

In this study we use a linguistic version of category theory in which the objects are drawn as text boxes describing some type of thing, like a protein or a genetic code, and where the arrows also have labels describing some functional relationship, as every protein has a genetic code. Chains of arrows can be composed, providing a description of how a number of small-scale relationships come together to constitute a single, conceptually simpler, larger-scale relation (like a person's father's sister's daughter is a simply their cousin; an example for "functionality" in the space of linguistics). These linguistic categories are called "ologs", short for "ontology logs" (see [47]). Ontology is the study of how or what something *is*, and ologs are a systematic framework in which to record the results of such as study. The term "log" (like a scientist's log book) alludes to the fact that such a study is never really complete, and that a study is only as valuable as it is connected into the network of human understanding. This brings us to the heart of the matter: in order to build a sufficient understanding of hierarchical materials, scientists must integrate their findings more precisely with those of other scientists.

The fact that an olog is fundamentally a category means that such connections can be formulated between ologs with mathematical rigor, and meaning preserved [47], to facilitate the communication with other fields of science. This concept is depicted for a simple example in Figure 2. Note how the structure of the category, *i.e.* the arrangement of objects (here: sets) and arrows (here: functions), is preserved while the objects and the arrows itself are subjected to a transformation. This means that if a certain property, such as the mechanical behavior of amyloids, can be described in a categorical framework, structure preserving transformations translate the components of the system into other systems, such as a wood or concrete based system, while the relations and thus the functionality within the category is maintained. The revelation and abstraction of the origin of protein material properties must be done by intensive materiomics studies that typically involve multiscale experiment and simulation.

We omit a precise definition of categories and hence ologs in this paper as we will focus on the application of this concept – the discussion will be limited to a general description of ologs and how they are constructed. Hence, we will proceed to describe ologs by example; for a more mathematically precise account of ologs, see reference [48].

**3. Results and discussion**

The use of ologs is a powerful tool that can ultimately enable the kinds of breakthroughs needed to further our understanding of how functional diversity is achieved despite intrinsic limitations of building blocks. The generation of ologs also allows us to observe the formation of patterns that define certain functionality, and draw connections between disparate fields. A key insight used here is that although patterns of functionality generation can be quite different in the space specific to applications (*e.g.*: proteins, language, music), they are remarkably similar in the space of categories. In other words, we hypothesize it is possible to observe universal patterns of how functionality is created in diverse fields; and that it is possible to generate fundamental laws (similar to PDEs in conventional physics) that describe the emergence of functionality from first principles.

We briefly expand on the potential powerful application of category theory mentioned above. As explained in Section 1 the same 20 amino acids can have different functions depending on how they are arranged in a sequence as defined by the genes. In other words, the same library of fundamental building blocks can produce different functionality depending on the precise sequence. Just so, an olog serves as a code or formula for a complex structure, but the context in which it is interpreted can lead to different results. We will show that the same olog can be interpreted as formulating the structural and functional relationship between a protein filament such as an alpha-helix and an amyloid fibril or the same relationships between two types of social networks involving a chain of participants. In the case of a protein the building blocks are polypeptide fragments or H-bond clusters as glue, whereas in a social network the building blocks are people and communication methods. It is the interplay between form and function of few universal building blocks that ties biological structuralism and category theory, and which may produce potentially novel approaches to designing engineered systems.

**3.1 Olog of two protein materials under axial loading: alpha-helix and amyloid fibril**

We develop an olog for two protein filaments that display a distinct mechanical behavior once exposed to mechanical force. We begin the discussion with a presentation of the proteins and their functional properties, here their mechanical properties under axial extension (realized *e.g. via* the application of an axial force applied to the protein filament). The structure, mechanisms and resulting functional properties have been developed in a series of earlier studies based on computational approaches to molecular nanomechanics (for alpha-helices, see [49] and for amyloids or beta-sheet crystals, see [20,50]; and we refer the reader to these original papers for further detail).

Figure 3 shows the visualization of the two protein materials based on an abstraction of how their mechanical properties can be understood based on the interplay of a set of "building blocks" (Figure 3A). Both protein materials resemble a linear arrangement of three elements, "bricks", "glue", and for one of them, "lifeline". Thereby as a design rule, brick and glue need to alternate in order to achieve a stable structure. Two brick or glue elements immediately next to each other would not stick together. There is a fundamental chemical reason for this constraint as bricks represent the protein's polypeptide backbone and glue represents H-bonding which can only occur between a cluster of amino acid residues in the backbone. The "lifeline" is a third element that is introduced here, reflecting the situation in which there is still a physical connection of bricks even after large force causes the glue to break. Chemically, this resembles the existence of a "hidden" polypeptide length such that there exists a "covalent" link between two brick elements even after the H-bond glue has broken. The hidden length is not observed as a relevant structural property until the glue breaks, at which point the lifeline comes into play and provides an increasing resistance against deformation. Thus, although both glue and lifeline can connect neighboring brick elements, they are differentiated in that the lifeline is much stronger than the glue and that its resting extension is roughly the failure extension of the glue (Figure 3).

Although this description of protein filaments is a simplification of how their mechanical properties can be described and the focus is set on a distinct feature of the protein material's behavior only, it enables us to understand the key functional properties based on the interplay of building blocks. We demonstrate this now with a detailed discussion of the two cases considered. Figure 3C depicts a model of an amyloid fibril (or similarly, model of a beta-sheet crystal as found in silk) subjected to axial deformation. The structure is realized by the assembly of on an alternating sequence of bricks (amino acid cluster) and glue (H-bond cluster). Upon the increase of the extension one of the glue elements breaks. Since there is no more physical connection between the two brick elements that were previously connected by the glue element the entire system has failed, and at an extension that is roughly equal to the failure extension of the glue (Figure 3E). We define this behavior as "brittle". Figure 3D depicts a model of an alpha-helix protein (a protein found for example in the cell's cytoskeleton) under axial loading, assembled based on an alternating sequence of bricks (amino acid cluster), glue (cluster of H-bonds) and a lifeline element. The lifeline element is formed by the protein backbone that still exists even after the cluster of H-bonds break after unfolding of one alpha-helical turn [49]; providing a physical connection that allows additional glue elements to break after more axial extension is applied. In fact, at large extensions all glue elements will have broken such that the system's overall failure extension is much larger than the failure extension of the glue, marking a "ductile" behavior (Figure 3E).

The comparison of the distinct mechanical behavior of alpha-helices and amyloid fibrils was achieved by mapping the key mechanisms that generate their specific properties into the abstract space of interactions between a set of building blocks. What was described in words in the preceding paragraphs can be rigorously achieved using ologs, which describe the interactions between building blocks. Through the development of ologs for each system we aim to answer a series of questions:

- What are the components of the system, and how do they interact?
- How do these interactions produce the functionality we observe of the overall system?
- When does functionality break down? *E.g.*, failure of building blocks as the system is pushed to extreme conditions, or the presence of defects.
- A further reaching question may be, by what process did the system come to be constructed, and what selective pressures at the macroscale induce observable changes in the system and at different levels in the structural makeup.

To eventually get us to this point, we will now discuss the components of our brick-and-glue system of proteins so as to acquaint the reader with the olog presented in Figure 4 which describes both the brittle and ductile protein filaments outlined above. Three universal elements, which we have been calling bricks (b), glue (g), and lifeline (L) are the abstract building blocks composing our systems, and they are defined in relation to one another as follows. Both glue and lifeline are materials that can connect two brick elements. There are two distinctions between them: *i*) the failure extension of glue is much less than that of brick, whereas the failure extension of lifeline is roughly equal to that of brick, and *ii*) the resting extension of lifeline is roughly equal to the failure extension of glue. These two properties ensure that the lifeline is not detected under axial loading until a glue element breaks and that all the glue elements break long before a lifeline or brick element breaks (see also Figure 3B).

This distinction between one number being *roughly equal* to another and one number being *much greater than* another is simple, yet universal in the sciences, and thus we can expect these types (M and O in the olog) to be quite common in scientific ologs. In fact, we reuse this concept within the olog when we distinguish a ductile system from a brittle one. That is, we characterize a ductile system to be one whose failure extension is *much greater than* that of its glue element, whereas we characterize a brittle system to be one whose failure extension is *roughly equal to* that of its glue element. Relations like that are typical for hierarchical systems where a scaling law applied to the scale of a building block connects the behavior of such building blocks to the overall system behavior. Other common (*i.e.* universal) patterns that we may find in biological materials is a certain shape (fibers, helices, spheres), bonds of a certain form (H-bonds, backbone), dimensionality (1D, 2D, 3D), and so on. Our olog concentrates on materials whose shape is one-dimensional, a feature we define by the use of mathematical graphs.

The interactions of building blocks are not limited by their interface. As each object represents a category itself, it can be again a composition of objects and arrows. Hence, the functionality can be affected by the existence or alteration of neighbor building blocks by drawing connections between objects within categories. This is exactly how functional and structural hierarchies are represented in an olog. Since a brick (or a glue) can refer to anything in the world, an entire system of bricks and glue can be regarded as a new "brick" (and a whole system of bonds as "glue"). A zoom-in or zoom-out is possible by defining new building blocks in terms of others.

Once the fundamental structure of our protein materials and the definition of ductility and brittleness have been defined in the olog, we describe our hypotheses by two arrows, 1:**A→E** and 5:**B→C**, the first of which hypothesizes that systems with lifelines are ductile, and the second of which hypothesizes that systems without lifelines are brittle. This hypothesis has now been examined in the paragraphs above, but can be even more carefully explicated using a category-theoretic formulation, where each component type and aspect is laid bare. In fact, we have no hope of proving an analogy between this protein setup and the upcoming social network setup without such a formulation. In Figure 4 we show the entire setup as a diagram of boxes and arrows, the precursor to an olog. However, this diagram is not sufficient in the sense that there are mathematical truths present in our system that are not present in the diagram. We include the rest of this information in Tables 1 and 2, which we will describe shortly.

In order to explain what is missing from Figure 4, we should more clearly explain what is there. Each box represents a set. For example box **H,** labeled "a graph", represents the set of graphs, whereas box **J**, labeled "a system consisting of bricks connected by glue, structured as in graph G", represents the set of such systems. Each arrow represents a function from one set to another, and its meaning is clear by reading the label of the source box, the label of the arrow, and then the label of the target box. For example, we read arrow 20: **J→H** as "a system consisting of bricks connected by glue, structured as in graph G is structured as a graph". Thus, any element of the set **J** is functionally assigned its structure graph, an element of **H**, by arrow 20. Just as the structure graph of a system is an *observable* of that system, any function from one set to another can be considered an observable of the former.

A function may be thought of as a "black box" which takes input of one type and returns output of another type. If the output of one function is fed as input to another function and the whole system is imbedded in a black box, it is called the composition of functions. Finally, two functions are equal (regardless of the inner workings of their "black boxes") if, upon giving the same input they always return the same output. The first kind of mathematical truth alluded to above that is missing from Figure 4 is a declaration of which compositions of functions in our system are equal. Such equalities of compositions of functions are called *commutative diagrams* in category theory literature. All such declarations are presented in Table 1. These equalities can be considered as checks on our understanding of all the sets and functions in the arrows – declaring them is at the very least "good science".

Table 2 describes a certain class of commutative diagrams; called *fiber product diagrams* (see also Figure S2). In a fiber product diagram, one set and two observables of it are declared as a kind of "universal solution" to a problem posed by another diagram. In these terms, we consider the diagram **D→H←J** as posing a problem, to which **D←F→J** is a solution, as we now explain. The diagram **D→H←J** poses the problem "what should we call a system consisting of bricks connected by glue, structured as in graph G, where graph G is a 'chain' graph". The declared solution is **F** "a one-dimensional system (S) of bricks (b) and glue (g)", together with its two observables **F→D** and **F→J**. Thus the second kind of mathematical truth alluded to above that is missing from Figure 4 is that some boxes and attributes have fixed meaning in terms of the others. A list of these is given in Table 2, where we see terms such as "one-dimensional", "brittle", "ductile", and "lifeline" defined solely in terms of more basic concepts.

Thus, while it is convenient to think of the olog for our protein systems as the diagram in Figure 4, in fact it is the totality of Figure 4, Table 1, and Table 2, which really constitute the olog. Just as in biological materials, the parts of the olog (its boxes and arrows) are not sufficient for the system to act as a whole – the less-obvious interrelationships between these parts give the system its functionality. It is important to note that ologs can be constructed based on modeling and simulation, experimental studies, or theoretical considerations that essentially result in the understanding necessary to formulate the olog. This has been done for the proteins considered here based on the results from earlier work which provided sufficient information to arrive at the formulation of the problem as shown in Figure 3.

**3.2 Olog of a social network**

In this section we construct a simple social network that may appear to some as vastly different as a protein filament, and to others as quite similar. The reason for the discrepancy is that semantically and physically the situations have almost nothing in common, but structurally and functionally they do. In fact, we will prove category-theoretically that they are structurally and functionally isomorphic in the sense that their ologs are identical. We now describe the setting for our simple social network as depicted in Figure 5. Imagine a building with sound-proof rooms labeled 1 through 100, equipped with a controlled wireless communication system connecting each pair of consecutive rooms. In each room a human participant sits on a chair with a simple wireless transceiver that can transmit and receive text messages from the participant to the left (his or her predecessor) or the person to the right (his or her successor). We assume that participants in odd numbered rooms are women and people in even numbered rooms are men, just for pronoun clarity. The goal is to faithfully pass messages (sentences of under ten words, say) from room 1 to room 100 and back the other way as quickly as possible. The woman in room 1 (respectively the man in room 100) receives a message from the experimenter. She then inputs it into her transceiver and sends it to her neighbor (2), who passes it along to his neighbor (3), and on down the line until it is received by the man in room 100, who submits it to the experimenter there. Thus the network has a task of faithfully sending messages from one experimenter to the other; if they fail to successfully transmit at least one message per hour we say that the system has failed.

An obstacle can be added by allowing that the transmission of messages between participants is not always error-free. That is, the experimenters can adjust the amount of "noise" in the system, resulting in messages that could be anywhere from error-free to completely unintelligible. For example, the message "the party was fun and exciting" may arrive in the next room as "tha partu was fon and escitin". Upon receiving a garbled message, a participant may take the time to "fix it up" before sending it along, thereby helping to ensure that the message can be correctly submitted at the end of the line. We define the "extension" of the system to be the amount of noise, measured as the probability that a transmission error occurs for an arbitrary letter in a message. Given sufficient noise, it may happen that no messages can be transferred successfully through the network. Thus, we define the "failure extension" of the network to be the amount of noise present when this occurs. Similarly, the failure extension of a glue element is the amount of noise at which a wireless transmission cannot be successfully sent from one room to the next.

Finally, we can add lifelines to this picture by adding physical passageways between consecutive rooms. Now, in case the noise gets too high, individuals may walk or run through these "lifeline passageways" and transmit a message by voice. During low levels of noise, the doorways will typically not be used to relay information because the text messaging is much faster, and hence the existence of the lifelines will be "hidden". However, once the transmission noise is severe enough to prevent good wireless communication (that is, the glue breaks), these passageways come into effect and save the network from breaking altogether. The three basic building blocks of this social network are shown in Figure 5A. For a rigorous analysis we also define a failure extension for bricks and lifeline, and resting extension for lifeline (Figure 5B). We define the failure extension of bricks and lifeline to be infinite (because messages existing on a given transceiver or passed *via* voice are unaffected by the noise level), and we define the resting extension of our lifeline passageways to be the amount of noise at which participants begin to use the passageways.

We now analyze the performance of the two types of networks constructed here, without and with a lifeline. In the system without a lifeline (Figure 5C) as soon as the noise level is high enough to cause breakdown of one of the glue elements the system fails since no more messages can be transmitted. In the system with lifelines (Figure 5D), even though glue elements may break there is still the possibility for signals to travel through the passageway such that a much greater noise level (or extension) can be sustained. A brittle network is one in which the failure extension is roughly the same as the failure extension for each glue element. A ductile network is one in which the failure extension is much greater than the failure extension of each glue element. We thus hypothesize that social networks with lifeline passageways will be ductile and that those without lifeline passageways will be brittle. While the above communication network is fairly degenerate as compared with, say the Facebook network, the basic idea is similar. People are connected with a set of "friends" and the basis of this friendship is communication. Communication can be muddled by various kinds of noise, but the use of additional forms of interaction (*e.g.* talking in face-to-face meetings in addition to using online text messages) increases the probability that the parties understand each other.

We have constructed a system so that the olog describing it is precisely the same as that defining the protein system of Section 2.1. The basic layout is in Figure 4, and Tables 1 and 2 add "rigidifying information". For example, the participants with their transceivers are the bricks, the wireless communication between neighboring rooms is the glue, the passageways are the lifelines. We define brittleness and ductility exactly as we did in the protein case and as described in the previous paragraph; in fact this is forced on us because boxes **C** and **E** are fiber products (see also Figure S2). The fact that the same olog describes our protein materials and our social network should be considered as a rigorous *analogy* or *isomorphism* between these two domains, as we describe in more detail in the next section.

### 3.3 Analogy between alpha-helical protein and communication network

The analogy between the protein strands (amyloid fibril and alpha helix) and the social network experiment is as follows: In both cases a network (protein/social) consisting of bricks (amino acid clusters/human participants) connected together by glue (H-bond cluster/wireless communication) is subject to pulling (axial extension/error-producing noise) and eventually reach a breaking point (when the maximum extension is reached/transmission rate drops to zero). Lifelines (additional physical connections *via* covalent links/passageways) serve to increase the ductility (failure extension of network divided by failure extension of individual glue elements/ditto) of the network. Table 3 gives a complete list of the meaning of all labels in the protein and social network.

We now rigorously show that the two situations can be modeled by precisely the same olog. Thus the olog sets out a space of possible systems that includes everything from proteins to social networks (and potentially many other realizations), any two instances of which must be analogous, at least to the level of detail found in Figure 4 and the associated tables. If one desires additional detail, for example to add a precise meaning for resting extensions, or even real numbers, one would simply expand the olog to capture these ideas. A key result from our discussion is that the interpretation of what *b, g* and *L* mean in different systems can be distinct (proteins, polymers, music, *etc.* can have different physical realizations of these concepts). Yet, their fundamental properties and how they relate to others – other elements, different scales in hierarchies, *etc.* – are defined properly in the olog, and mathematically expressed not only as a fundamental property but in addition as functors to other elements in the system. For physical systems a key aspect of understanding the interplay of building blocks can for instance be expressed in scaling laws that define properties as a function of ratios of length-scales or energy levels, which fundamentally define how elements behave and interact with others. The general presentation of such relationships in networks is what is missing in current theories, and is where ologs present a powerful paradigm for *de novo* design of biologically inspired systems that span multiple hierarchical levels. This is because ologs achieve a rigorous description of the synergistic interactions of structures and mechanisms at multiple scales, which provides the basis for enhanced functionality despite the reliance on few distinct building blocks.

It is important to note that for the sake of the analogy discussed here the two very different domains (protein *vs.* social network) were *designed* in a way to show that they could have identical conceptual descriptions at a very high level of detail. More detail could show differences between these two domains. For example, an observation we purposely did not include is that the bricks in our social network need to breathe and eat. It is impossible (and perhaps not desirable in some cases) to include every detail of each system – our goal was to emphasize the essential parameters, and to provide a level of abstraction that emphasizes the key elements that define functional properties. Furthermore, whereas it may be rare for two different scientific studies to result in identical ologs, finding reusable parts should be quite common. In our olog, the notion of bricks being connected together by glue to form the structure of a graph is surely reusable not only within materials science but throughout engineering.

Of course, the biological system was strongly simplified and we focused only on a single aspect within the vast range of properties in biological materials. Similarly, our social network was contrived to fit the olog of the protein. This analysis does not claim to describe the formation of "all" complex protein assemblies but shall serve as a generally comprehensible example for a methodology which puts these kinds of analogies in a concise mathematical framework, where future work could emphasize on applications to more complex cases. It shall thus serve as an impetus for further studies in this field. In fact, without the classification as a category theoretic transformation, such analogies have been recently used to compare active centers of proteins, that means a cluster of amino acids that have a high centrality in the amino acid network of the hosting protein by their participation in enzyme catalysis and substrate binding, to "strangers" in social networks and top predators in mammalian networks [51]. Creative elements, a highly specific subset of these central residues, occupy a central position in protein structure networks. They, among other things, give non-redundant, unique connections in their neighborhood, integrate the communication of the entire network, and accommodate most of the energy of the whole network. In atypical situations they become especially relevant due to their transient, weak links to important positions such as hubs. In mammalian networks top predators take the role of active centers as they act as couplers of distinct and dissimilar energy channels and increase the stability of the ecosystem's network whereas in social networks "strangers" – often innovators and successful managers – occupy "structural holes". They show exactly the same functional behavior as the active centers in protein structure networks and thus these networks are connected by structure preserving transformations or *functors*.

While we cannot discuss it here in detail (as it would be out of scope of this article), the category-theoretic notion of *functors*, which formally connect one olog to another, will eventually allow scientists in vastly different fields to share their work by rigorously connecting together their ologs. This opens enormous opportunities for design of novel functional properties by drawing from the understanding gained in diverse fields. Observations made in one field, *e.g.* about the dynamic response and transformation of active centers after repeated stress and the re-organization of the network topology, lead to insights in other fields, *e.g.* about the flow of novel information and direction of evolution in companies with well-connected collaborators whose contact structure increases the performance in uncertain environments or in crisis [51,52,53]. Other recent studies of the characteristics of biological protein networks showed that the modular structure of these networks improves the robustness against hub malfunction while increasing vulnerability against random failure which stands in contrast to the behavior of typical non-biological networks such as the Internet [54]. Insights gained from these types of studies accompanied by a systematic description of its functional features may help in the construction of artificial networks that inherit the advantageous properties of the biological archetype.

**4. Conclusion**

A unique aspect of the analysis provided here is that we described a rigorous analysis of the conceptual interaction rules in protein materials and establish a direct link to those of a social network *via* the use of category theory. This qualitative account will allow us to draw direct analogies to existing models of complex hierarchical structures such as those from social networks, and potentially linguistic theory where similar problems have been studied, and enables the utilization of insights and design paradigms across disparate fields of the science of hierarchical systems (Figure 6). The key concepts presented provides a generic framework that has the potential to unify existing understanding derived from the myriad of existing studies of individual protein materials such as bone, silk, or cells and many others, where a major limitation was that no unifying framework that applies generally to all such materials has yet been proposed. This paradigm and associated design rules, which are applicable to other complex systems such as music, engineered technology and materials, or food recipes, could emerge as an exciting new field of study and make critical contribution to the field of materiomics for which it serves as a central tool to describe structure-function relationships for hierarchical systems.

Future directions, open research questions, and the impact of an increased understanding of hierarchical protein materials is discussed at three levels with increasing generality: *i*), impact on protein material synthesis (design, engineering and manufacturing or novel biomaterials), *ii*), impact on bioinspired nanoscale material design and assembly (*e.g.* hierarchical materials such as fibers, yarns or armors), and *iii*), impact on macro-scale systems design and engineering (*e.g.* design of cars, airplanes *etc.* where the merger of the concepts of structure and material across all the scales provides opportunities for more efficient systems). Immediate future work could be directed towards applying the concept of ologs to specific hierarchical biological materials, such as to silk or bone that show a greater complexity than the simple problems reviewed here. While the resulting ologs are more complex, the basic approach is identical and the main insights discussed here should still hold. Eventually, the olog shown in this paper (Figure 4) could be implemented in a computational model, which will open the possibility for design optimization using numerical algorithms or make it easier to reuse existing ologs for the design of new ones.

**Acknowledgements**

MJB acknowledges support from a DOD-PECASE award. DIS acknowledges support from Office of Naval Research grant (N000141010841) as well as generous support from Clark Barwick, Jacob Lurie, and the Massachusetts Institute of Technology Department of Mathematics. TG acknowledges support from the German Academic Foundation (Studienstiftung des deutschen Volkes). We acknowledge helpful discussions with T. Ackbarow. The funders had no role in study design, data collection and analysis, decision to publish, or preparation of the manuscript. The authors declare that no competing interests exist.

**Supplementary material**

Supplementary material is available on the journal website.

**Figures & figure legends**

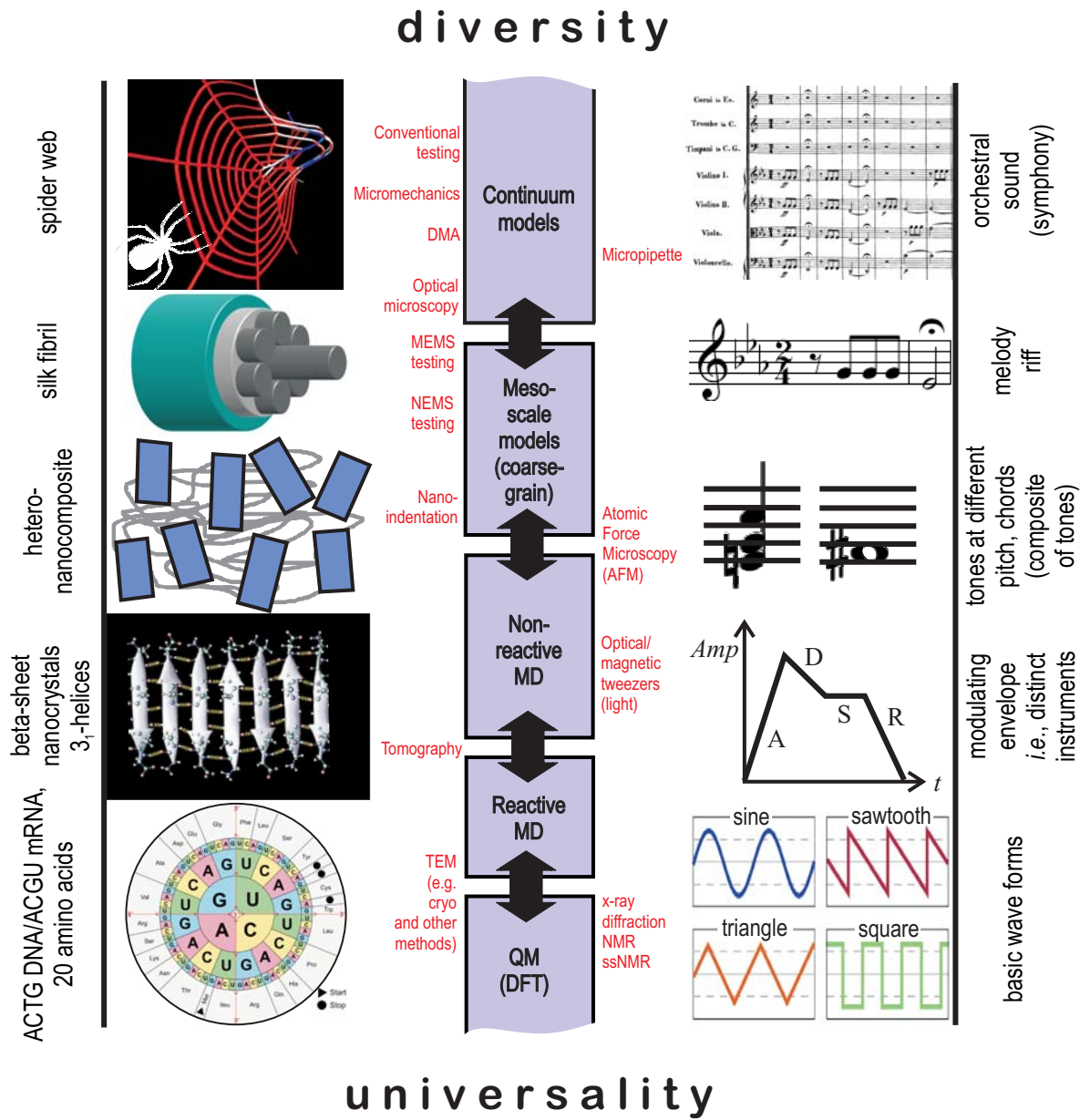

Figure 1. **Illustration of multiscale hierarchical structure of protein materials, a summary of multiscale modeling and experimental tools, and an analogy to music** (figure adapted from [14]). In protein materials (left for the example of spider silk), multifunctional materials are created *via* the formation of hierarchical structures. The synergistic interaction of structures and mechanisms at multiple scales provides the basis for enhanced functionality of biological materials despite the reliance on few distinct building blocks. Similar to the case of protein materials is music (right), where universal elements such as basic wave forms or a set of available instruments are used in hierarchical assemblies to provide macroscale functionality, and eventually a particular orchestral sound (*e.g.* a symphony). Universality tends to dominate at smaller levels, whereas diversity is found predominantly at larger, functional levels. The integrated use of computational and experimental methods at multiple scales provides a powerful approach to elucidate the scientific concepts underlying the materiomics paradigm (center).

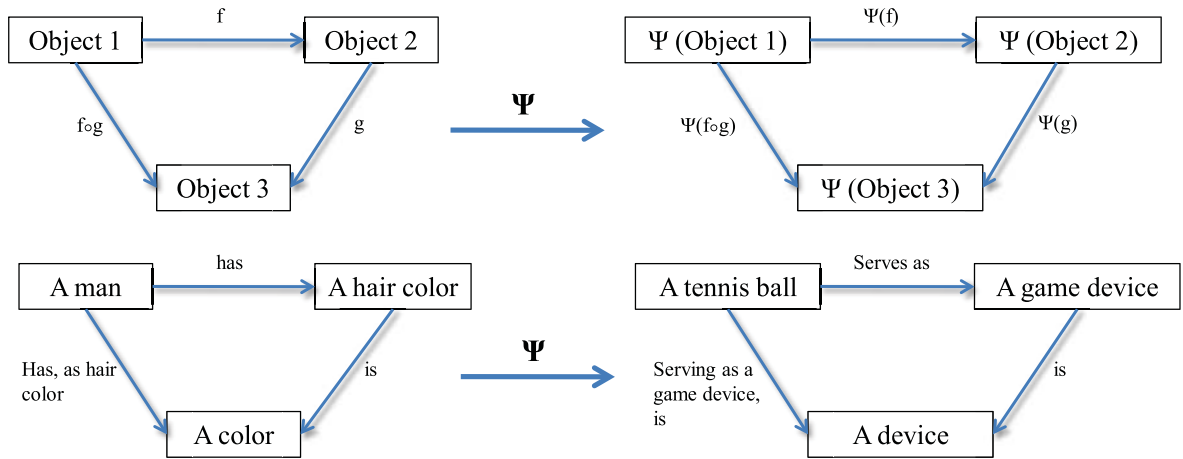

**Figure 2. Simple examples of transformations preserving structure in category theory**. Categories consist of objects and arrows which are closed under composition and satisfy certain conditions typical of functions. Ψ is a structure-preserving transformation (or covariant *functor*, or *morphism* of categories) between the two categories. If the categorical objects in this example are considered as sets of instances, then each instance of the set 'A man' is mapped to an instance of the set 'A tennis ball'. This concept applies to all objects and arrows in the categories.

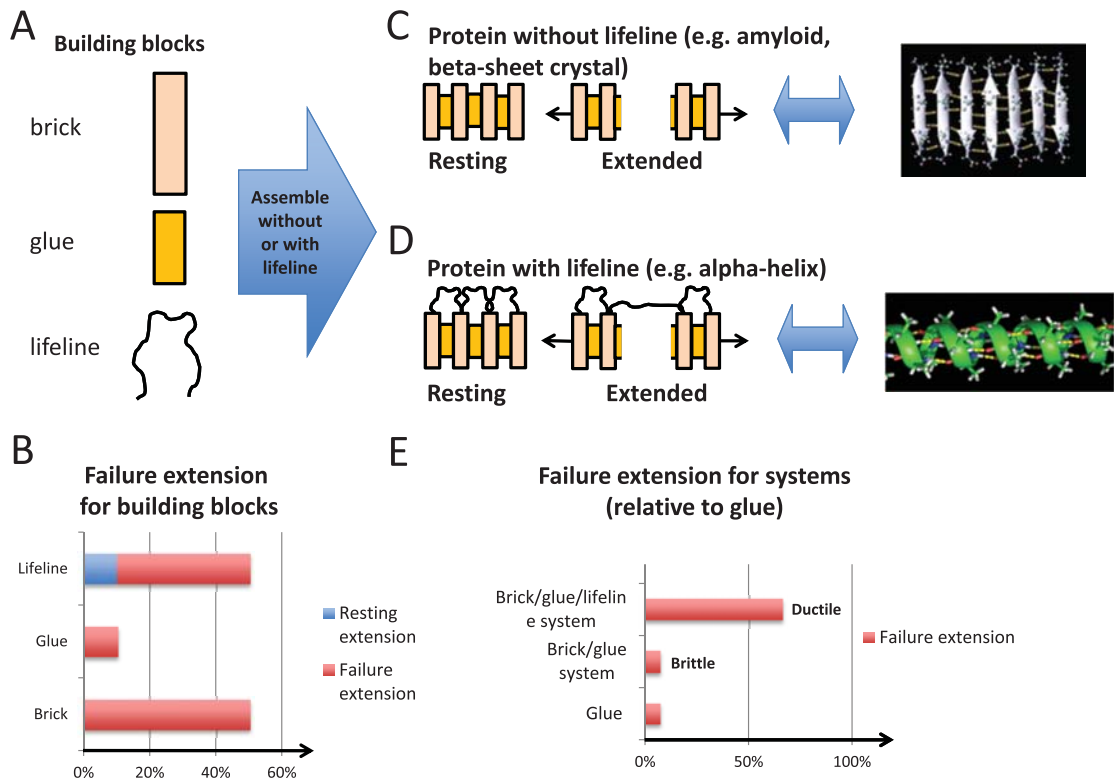

**Figure 3. Visualization of protein materials considered here, and abstraction of how key functional properties (here: mechanical properties under axial extension) can be understood based on the interplay of a set of "building blocks".** A, Overview over fundamental building blocks of our protein materials. The protein materials considered here are composed of a linear arrangement of three elements, "bricks", "glue", and in some cases "lifeline". Thereby as a "design rule", brick and glue need to alternate in order to achieve a stable structure. That is, two brick or glue elements immediately next to each other would not stick together – the chemical reason is that bricks represent the protein's polypeptide backbone and glue represents H-bonding which can only occur between residues in the backbone. The "lifeline" is a third element introduced here, reflecting the situation when there is still a physical connection between bricks even after the glue breaks. Chemically, this resembles the existence of "hidden" polypeptide length such that there exists a "covalent" link between two brick elements even after the H-bond glue has broken, where the hidden length is not "visible" before the glue is actually broken. B, Mechanical behavior of each of the building blocks characterized by a description of the failure extension. The hidden length of lifelines is reflected in the fact that the resting extension of the lifeline is roughly equal to the failure extension of the glue. Both the brick and the lifeline have large failure extensions relative to the glue. C, Model of an amyloid-like protein (or similarly, a beta-sheet crystal as depicted in the image on the right) under axial loading. This resembles a system without a lifeline since after breaking of the H-bond cluster (= glue) between the layers formed by clusters of polypeptide (= brick) no physical connection exists. D, model of an alpha-helical protein under axial loading. This resembles a system with a lifeline, as after breaking of the cluster of ≈3-4 H-bonds (= glue) that are formed between clusters of amino acids (=brick) there still exists a physical connection due to the polypeptide backbone as shown in D (= lifeline). As shown in E, the existence of a lifeline has major implications on the functional properties of the overall system. A system with a lifeline (D) shows a ductile response, where a connection can be sustained at large extension as compared to the glue alone. In contrast a system without a lifeline (C) shows a brittle response, where only a small extension can be sustained until the material breaks (which equals roughly the failure extension of the glue).

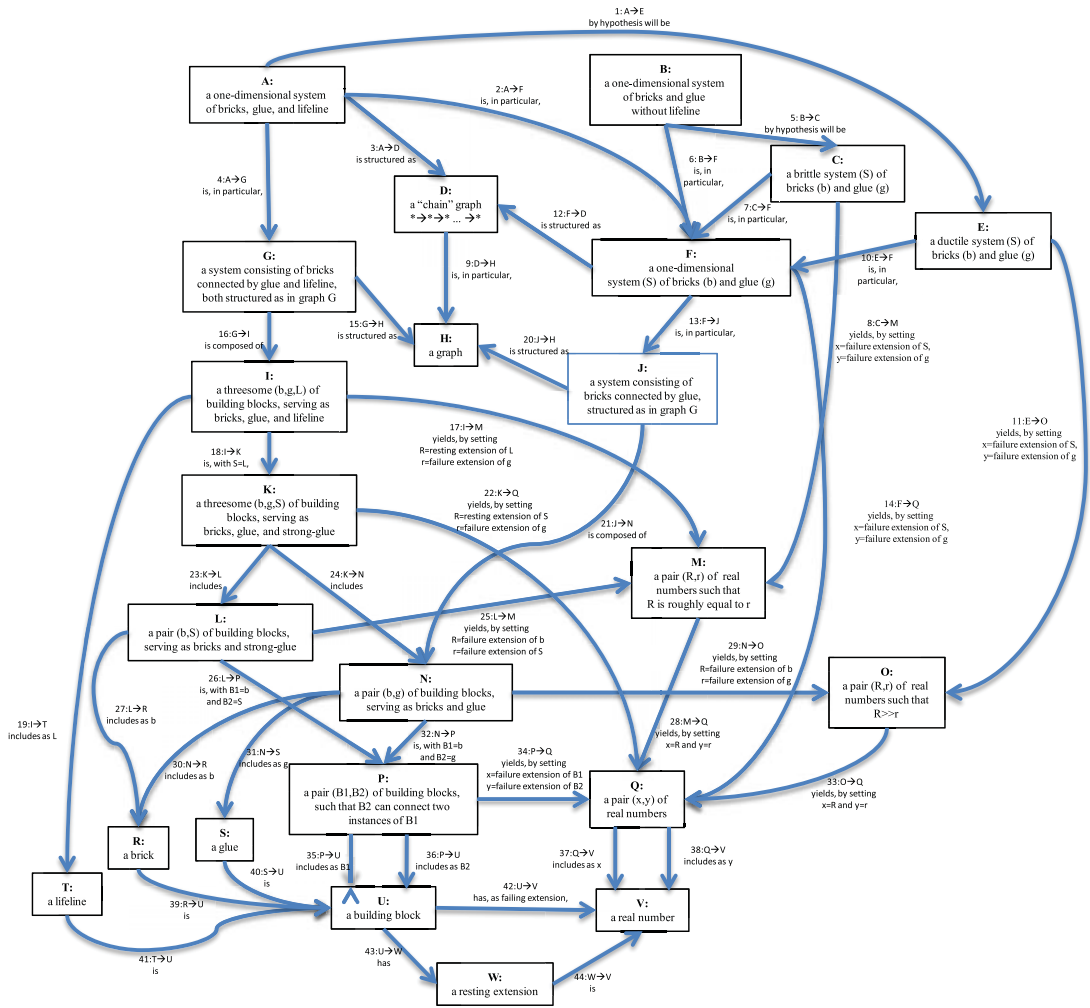

**Figure 4**. **Pictured here is an olog, which captures the semantic content of our situations, as described in Sections 3.1 – 3.3**. Each box represents an abstract type, and each arrow represents an aspect (or observable) of that type. Each type refers to a set of intended instances, which we think of as being contained in the box. For example, box **E** contains ductile sequences of bricks and glue (like an alpha helix), whereas box **V** contains real numbers (like 9.228). Each arrow from a source box to a target box refers to an observation one may make on things in the source box, for which the observed result is a thing in the target box. For example, arrow 11:**E**→**O** indicates that one can observe of any ductile material S a pair of numbers (R,r) where R is much greater than r. The meaning of these numbers R and r is enforced by a "commutative diagram" declared in Table 1 (line 6): the number R must refer to the failure extension of the system S and the number r must refer to the failure extension of its glue. This says that a ductile system fails at a much greater extension than its glue elements do. Perhaps a simpler but more mundane observation is made by arrow 37:**Q**→**V** which indicates that one can take any pair of real numbers (*x,y*) and observe the *x*-coordinate. So the pair (8.0, 3.2) is inside box **Q**, and it is observed by arrow 37 to have *x*-coordinate 8.0, which is in box **V**. Thus, each box is meant to contain an intended set of instances and each arrow is meant to functionally relate two such sets. The rest of the olog is recorded in Tables 1 and 2. Some are commutative diagrams which declare two paths through the olog to be equivalent and some are fiber products which define new types in terms of others.

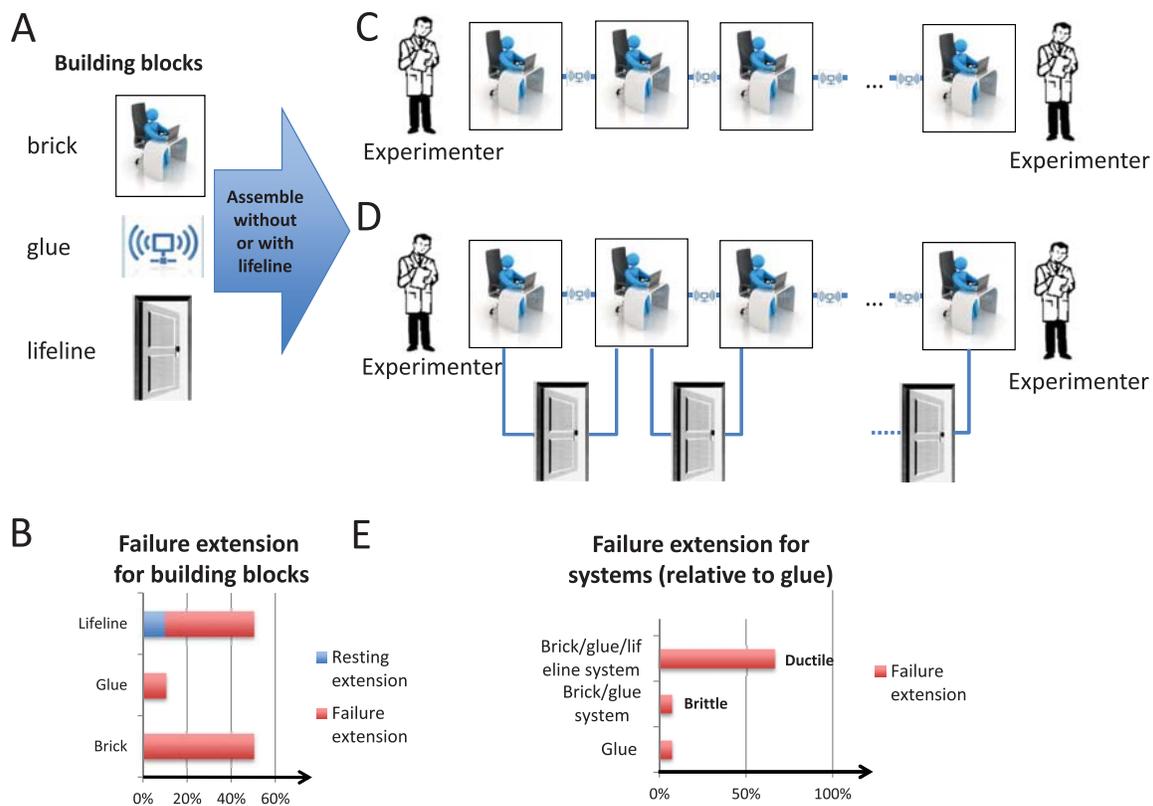

**Figure 5. Visual representation of the social network**. A, Overview over fundamental building blocks of our social networks. The social networks considered here are composed of a linear arrangement of three elements, "bricks", "glue" and in some cases, "lifeline". Thereby, as a "design rule", brick and glue need to alternate in order to achieve a stable structure. That is, two brick or glue elements immediately next to each other would not stick together; where the chemical reason is that bricks represent participants with transceivers and glue represents wireless communication that, in our case, can only occur between neighboring participants. The "lifeline" is a third element that is introduced here, reflecting the situation when there is still a physical connection of bricks even after the glue breaks. This reflects the existence of a "hidden" connection in that there exists a physical passageway between two brick elements even after the communication over the wireless connection is no longer feasible. The hidden connection is not "visible" before the glue is actually broken because, for reasons of efficiency, participants will choose to communicate the simple messages wirelessly rather than verbally, as the latter requires much more effort. B, Mechanical behavior of each of the building blocks. The hidden length of lifelines is reflected in the fact that the resting extension of the lifeline is roughly equal to the failure extension of the glue. In other words, lifeline passageways are used only when wireless communication is no longer feasible. Both the brick and the lifeline have large failure extensions relative to the glue because participants and their verbal communication function perfectly well in the presence of noise on the wireless channels. C, Representation of a *social network not allowing for face-to-face interaction* under stress from wireless noise. This resembles a system without a lifeline, as after noise on the wireless line reaches a critical point, messages can no longer be correctly conveyed. D, Representation of a *social network allowing for face-to-face interaction* under high levels of wireless noise. This resembles a system with a lifeline, as after messages can no longer be conveyed wirelessly, communication can still take place, due to the physical passageways as shown in D. As shown in E, the existence of a lifeline has major implications on the functional properties of the system. A system with a lifeline (D) shows a ductile response, where a connection can be sustained at large displacements as compared to the glue alone. In contrast a system without a lifeline (C) shows a brittle response, where only a small displacement can be sustained until the material breaks (roughly the failure extension of the glue).

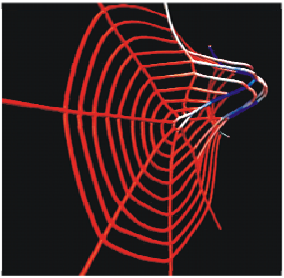 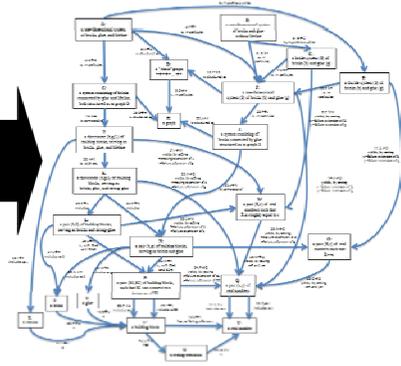 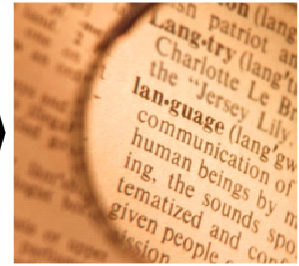

**Figure 6. Schematic illustration of the approach discussed here, the representation of complex hierarchical systems such as biological materials (*e.g.* silk) and language in the same category theory space (olog).** The description of how functional properties emerge in different hierarchical systems can be rigorously described using this approach, and fundamental insight can be derived. This finds immediate applications in the design of synthetic systems (*e.g.* novel materials).

**Tables**

| Starting point | Ending point | Path 1 | Path 2 | Same result |
|---|---|---|---|---|
| A: a one-dimensional system of bricks, glue, and lifeline | F: a one-dimensional system (II) of bricks (b) and glue (g) | A →1→ E →10→ F | A →2→ F | Each of these paths from A to F simply "forgets" the lifeline. |
| A: a one-dimensional system of bricks, glue, and lifeline | D: a "chain" graph •→•→•→ →• | A →2→ F →12→ D | A → D | Each of these paths from A to D yields the structure graph of the system, which is a "chain" graph. |
| A: a one-dimensional system of bricks, glue, and lifeline | H: a graph | A →3→ D →9→ H | A →4→ G →15→ H | Each of these paths from A to H yields the structure graph of the system. |
| B: a one-dimensional system of bricks and glue without lifeline | F: a one-dimensional system (II) of bricks (b) and glue (g) | B →5→ C →7→ F | B →6→ F | Each of these paths from B to F simply forgets that the system has no lifeline. |
| C: a brittle system (II) of bricks (b) and glue (g) | Q: a pair (x,y) of real numbers | C →7→ F →14→ Q | C →8→ M →28→ Q | Each of these paths from C to Q sets x=failure extension of the system (II), y=failure extension of the glue (g). |
| E: a ductile system (II) of bricks (b) and glue (g) | Q: a pair (x,y) of real numbers | E →10→ F →14→ Q | E →11→ O →33→ Q | Each of these paths from E to Q sets x=failure extension of the system (II), y=failure extension of the glue (g). |
| F: a one-dimensional system (II) of bricks (b) and glue (g) | H: a graph | F →12→ D →9→ H | F →13→ J →20→ H | Each of these paths from F to H yields the structure graph of the system. |
| I: a threesome (b,g,L) of building blocks, serving as bricks, glue, and lifeline | Q: a pair (x,y) of real numbers | I →17→ M →28→ Q | I →18→ K →22→ Q | Each of these paths from I to Q sets x=resting extension of lifeline (L), y=failure extension of glue (g). |
| I: a threesome (b,g,L) of building blocks, serving as bricks, glue, and lifeline | U: a building block | I →19→ K →23→ L →26→ P →36→ U | I →19→ T →41→ U | Each of these paths from I to U yields the lifeline element (L). |
| K: a threesome (b,g,H) of building blocks, serving as bricks, glue, and strong-glue | R: a brick | K →24→ N →30→ R | K →25→ L →27→ R | Each of these paths from K to R yields the same brick element (b). |
| L: a pair (b,H) of building blocks, serving as bricks and strong-glue | Q: a pair (x,y) of real numbers | L →26→ P →34→ Q | L →25→ M →28→ Q | Each of these paths from L to Q sets x=failure extension of brick (b), y=failure extension of strong-glue (H). |
| L: a pair (b,H) of building blocks, serving as bricks and strong-glue | U: a building block | L →26→ P →35→ U | L →27→ R →39→ U | Each of these paths from L to U yields the brick element (b). |
| N: a pair (b,g) of building blocks, serving as bricks and glue | Q: a pair (x,y) of real numbers | N →29→ O →33→ Q | N →32→ P →34→ Q | Each of these paths from N to Q sets x=failure extension of brick (b), y=failure extension of glue (g). |
| N: a pair (b,g) of building blocks, serving as bricks and glue | U: a building block | N →32→ P →36→ U | N →31→ S →40→ U | Each of these paths from N to U yields the brick element (b). |
| N: a pair (b,g) of building blocks, serving as bricks and glue | U: a building block | P →34→ Q →37→ V | P →35→ U →42→ V | Each of these paths from N to U yields the glue element (g). |
| P: a pair (B1,B2) of building blocks, such that B2 can connect two instances of B1 | V: a real number | P →34→ Q →38→ V | P →36→ U →42→ V | Each of these paths from P to V yields the failure extension of B1. |
| P: a pair (B1,B2) of building blocks, such that B2 can connect two instances of B1 | V: a real number | | | Each of these paths from P to V yields the failure extension of B2. |

**Table 1. Commutative diagrams in the olog.** Each sequence of consecutive arrows through the olog (Figure 4) is called a path, which represents a functional relationship between its starting point and its ending point. Two such paths A→B may result in the same function, and the 17 lines of this table record 17 cases of this phenomenon in our olog. The idea is that given an instance of A, each of these paths returns the same instance of type **B**. By having this additional data, we confine the meaning of the label on each box and arrow – they cannot stray far from our intended meaning without "breaking" these path equalities. Thus this table serves as an additional check on our labels. [For a more diagrammatic description of the same information presented in the typical style of category theory, see Figure S1.]

| Object | Fiber product object name | Defining attributes | Equated terms | "Idea" |
|---|---|---|---|---|
| A | a one-dimensional system of bricks, glue, and lifeline | D →3 A →4 G | D →9 H ←15 G | A system of bricks, glue, and lifeline is defined as "one-dimensional" if its structure graphs (brick/glue) and (brick/lifeline) are both chains. |
| C | a brittle system of bricks (b) and glue (g) | F ←7 C →8 M | F →14 Q ←28 M | A system is defined as "brittle" if its failure extension is roughly equal to the failure extension of its glue. |
| E | a ductile system of bricks (b) and glue (g) | F ←10 E →11 O | F →14 Q ←33 O | A system is defined as "ductile" if its failure extension is much greater than the failure extension of its glue. |
| F | a one-dimensional sequence (S) of bricks (b) and glue (g) | D →12 F →15 J | D →9 H ←20 J | A system of bricks and glue is defined as "one-dimensional" if its structure graph is a chain. |
| I | a threesome (b,g,L) of building blocks, serving as bricks, glue, and lifeline | M ←17 I →18 K | M →28 Q ←27 K | A strong-glue element is defined as "lifeline" if its resting extension is roughly equal to the failure extension of a glue element. |
| K | a threesome (b,g,S) of building blocks, serving as bricks, glue, and strong-glue | N ←24 K →23 L | N →30 R ←27 L | A "brick/glue/strong-glue threesome" is defined to be a brick/glue pair and a brick/lifeline pair where the bricks are the same in both instances. |
| L | a pair (b,S) of building blocks, serving as bricks and strong-glue | M ←25 L →26 P | M →28 Q ←34 P | Two building blocks, one of which can connect together two instances of the other, are defined as "bricks and strong-glue" if their failure extensions are roughly equal. |
| N | a pair (b,g) of building blocks, serving as bricks and glue | O ←29 N →32 P | O →33 Q ←34 P | Two building blocks, one of which can connect together two instances of the other, are defined as "bricks and glue" if the failure extension of the connector is much less than the failure extension of the connectee. |

**Table 2. Fiber product diagrams in the olog.** Some boxes in the olog (Figure 4) are defined in terms of others by use of so-called fiber products. For example, object **A** is defined in terms of three others in relationship, **D→H←G**: given a system of bricks, glue, and lifeline (**D**), we observe its structure graph (**H**) and set it equal to a "chain" graph (**G**) – in so doing we define "one-dimensionality" for a system. A reader of this olog realizes that our notion of one-dimensionality is not up for interpretation, but directly dependent on the other notions in this olog. By having this additional data, we confine the meaning of 24 labels (8 for boxes, 16 for arrows) in the olog. Thus this table serves to anchor the interpretation of our olog more firmly to its original intention. [For a more diagrammatic description of the same information presented in the typical style of category theory, see Figure S2.]

| Type | Type Labels | Protein Specific | Social-network Specific |
|---|---|---|---|
| A | a one-dimensional system of bricks, glue, and lifeline | alpha-helix | social network with wireless & physical passageways |
| B | a one-dimensional system of bricks and glue without lifeline | amyloid | social network with wireless, without physical passageways |
| C | a brittle system (S) of bricks (b) and glue (g) | brittle protein filament | brittle social network |
| D | a "chain" graph *→*→* ... →* | chain shape for protein | chain shape for network |
| E | a ductile system (S) of bricks (b) and glue (g) | ductile protein filament | ductile social network |
| F | a one-dimensional system (S) of bricks (b) and glue (g) | alpha-helix / amyloid | social network |
| G | a system consisting of bricks connected by glue and lifeline, both structured as in graph G | lifeline protein of specified shape | lifeline social network of specified shape |
| H | a graph | shape of protein | shape of network |
| I | a threesome (b,g,L) of building blocks, serving as bricks, glue, and lifeline | amino cluster, H-bond, backbone | transceiver, wifi system, physical passageway |
| J | a system consisting of bricks connected by glue, structured as in graph G | protein of specified shape | social network of specified shape |
| K | a threesome (b,g,S) of building blocks, serving as bricks, glue, and strong-glue | amino acid cluster, H-bond, backbone | transceiver, wifi system, physical passageway |
| L | a pair (b,S) of building blocks, serving as bricks and strong-glue | amino acid cluster, backbone | transceiver, physical passageway |
| M | a pair (R,r) of real numbers such that R is roughly equal to r | e.g. R=20.5 r=23 .45 | e.g. R=20.5 r=23 .45 |
| N | a pair (b,g) of building blocks, serving as bricks and glue | amino acid cluster, H-bond | transceiver, wifi system |
| O | a pair (R,r) of real numbers such that R>>r | e.g. R=100 r=20.6 | e.g. R=100 r=20.6 |
| P | a pair (B1,B2) of building blocks, such that B2 can connect two instances of B1 | e.g. amino acid and backbone | e.g transceiver and wifi |
| Q | a pair (x,y) of real numbers | e.g. x=20.55, y=50.6 | e.g. x=20.55, y=50.6 |
| R | a brick | amino acid cluster | transceiver |
| S | a glue | H-bond cluster | wifi connection |
| T | a lifeline | backbone | physical passageway |
| U | a building block | basic unit of material | basic unit of social interaction |
| V | a real number | e.g. 181.2 | e.g. 181.2 |
| W | a resting extension | e.g. 61 Angstrom | e.g. 1/100 error/bit |

**Table 3. Analogy between protein and social network.** Because our olog (Figure 4) was designed to abstract away the particulars of either the protein or the social network (using terms like "brick" instead of "amino-acid cluster" or "transceiver"), this table serves to remind the reader of the particulars in each case. Each type in the olog is described in these two cases. Some types, such as "a real number", stand on their own and we merely give examples. Others, such as "a one-dimensional system of bricks, glue, and lifeline" require a bit more description in the concrete cases. For more on this, see sections 2.1 through 2.3. This table provides the necessary description to connect the concrete formulations in the case of our protein and social network to the abstract formulation given by Figure 4.

**Supplementary Figure**

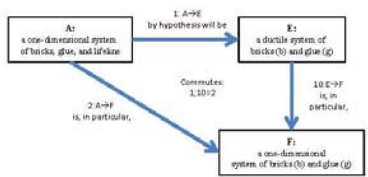 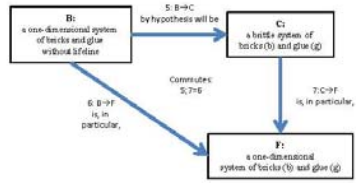

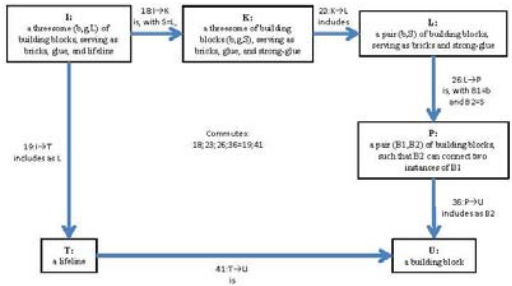 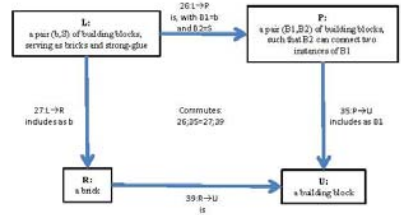

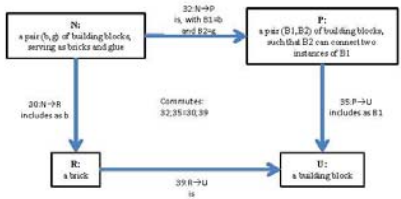 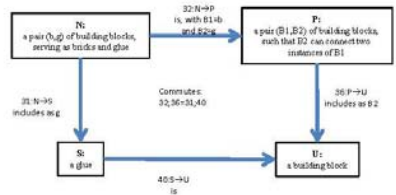

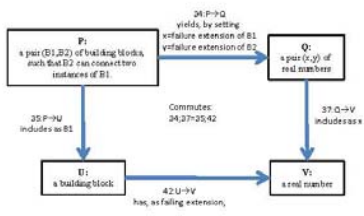 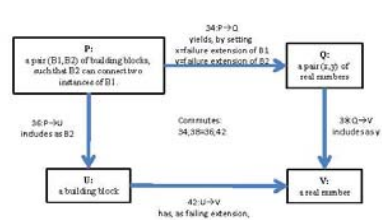

**Figure S1. Commutativity in the olog of the protein**. In each of these eight diagrams, there are two paths from the upper left-hand box to the lower right-hand box. By stating that these diagrams are commutative, we are saying that these two paths are equivalent – given the same input they produce the same output. For example it is declared 30;39=32;35 : **N→U**, which means that starting with a pair (b,g) of building blocks serving as bricks and glue, one can obtain a building block in two ways, but either way the answer is the same: the brick. Similarly, 31;40=32;36 : **N→U**, which means that again starting with (b,g) we can again obtain a building block in two ways, but either way the answer will be the same: glue. An example of a non-commutative diagram found in the original olog is: 31;40 ≠ 30;39 : **N→U**. Starting with a pair (b,g), the path 31;40 produces its glue element whereas the path 30;39 produces its brick element. These facts are in some sense obvious, but to make ologs a rigorous system such facts must be recorded.

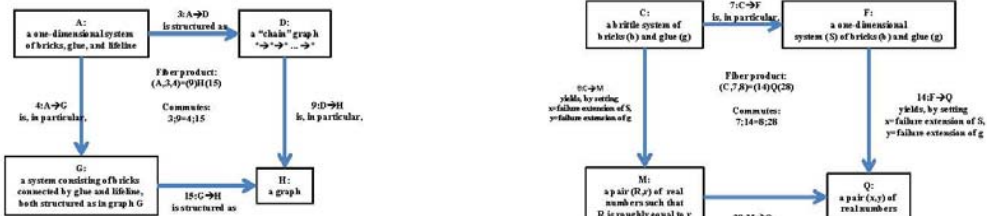
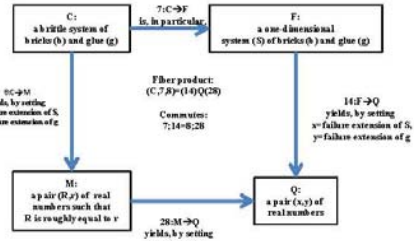
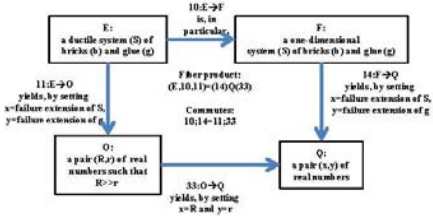
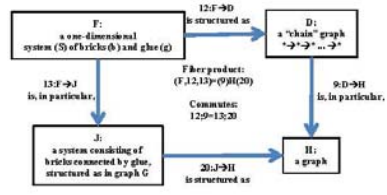
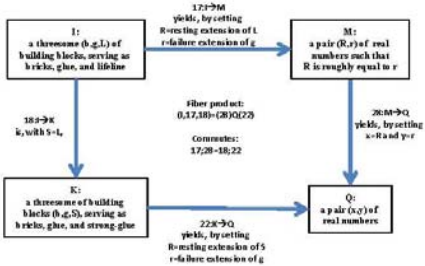
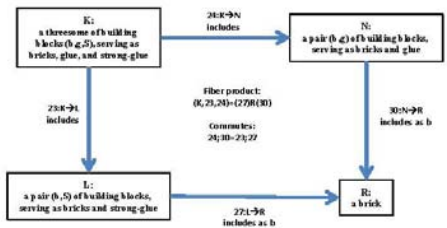
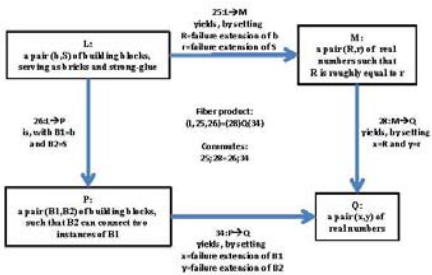
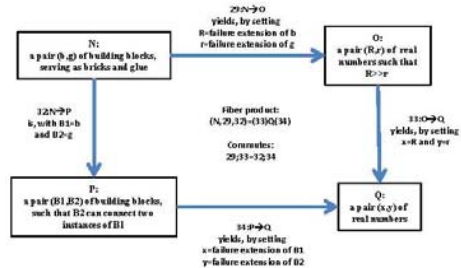

**Figure S2. Fiber products in the olog of the protein.** In each case, the upper left-hand box is the "fiber product" of the rest of the square. The property of being a fiber product defines the upper left-hand object: for example the notion of "one-dimensionality" in box **A** is defined for a system of bricks, glue, and lifeline by examining the structure of that system as a graph, and forcing that this graph is a chain graph (i.e. the elements are connected one to the next in a line).

23